\documentstyle[12pt]{article}
\begin{document}
\newpage

\setcounter{page}{0}
%
%
%
%
\def\CC{{\Bbb C}}
\def\NN{{\Bbb N}}
\def\QQ{{\Bbb Q}}
\def\RR{{\Bbb R}}
\def\ZZ{{\Bbb Z}}
\def\cA{{\cal A}}          \def\cB{{\cal B}}          \def\cC{{\cal C}}
\def\cD{{\cal D}}          \def\cE{{\cal E}}          \def\cF{{\cal F}}
\def\cG{{\cal G}}          \def\cH{{\cal H}}          \def\cI{{\cal I}}
\def\cJ{{\cal J}}          \def\cK{{\cal K}}          \def\cL{{\cal L}} 
\def\cM{{\cal M}}          \def\cN{{\cal N}}          \def\cO{{\cal O}}
\def\cP{{\cal P}}          \def\cQ{{\cal Q}}          \def\cR{{\cal R}} 
\def\cS{{\cal S}}          \def\cT{{\cal T}}          \def\cU{{\cal U}}
\def\cV{{\cal V}}          \def\cW{{\cal W}}          \def\cX{{\cal X}}
\def\cY{{\cal Y}}          \def\cZ{{\cal Z}}
\def\qed{\hfill \rule{5pt}{5pt}}
\newtheorem{lemma}{Lemma}
\newtheorem{prop}{Proposition}
\newtheorem{theo}{Theorem}
\newenvironment{result}{\vspace{.2cm} \em}{\vspace{.2cm}}

\rightline{CPTH-9611477}  
\rightline{q-alg/9611015} 
\rightline{November 96} 

\vfill
\vfill

\begin{center}

   {\LARGE   {\bf {\sf 
A Generalization of ${\cal U}_{h}(sl(2))$ via Jacobian Elliptic
Functions}}}\\[2cm]

\smallskip 

{\large A. Chakrabarti\footnote{chakra@orphee.polytechnique.fr}}

\smallskip 

\smallskip 

\smallskip

{\em  \footnote{Laboratoire Propre du CNRS UPR A.0014}Centre de Physique 
Th\'eorique, Ecole Polytechnique, \\
91128 Palaiseau Cedex, France.}

\end{center}

\vfill

\begin{abstract}
A two-parametric generalization of the Jordanian deformation 
${\cal U}_{h}(sl(2))$
of $sl(2)$ is presented. This involves Jacobian elliptic functions. In our
deformation ${\cal U}_{(h,k)}(sl(2))$, for $k^2=1$ one gets back ${\cal U}_{h}(sl(2))$.
The constuction is presented via a nonlinear map on $sl(2)$. This invertible
map directly furnishes the highest weight irreducible representations of
${\cal U}_{(h,k)}(sl(2))$. This map also provides two distinct induced Hopf stuctures,
which are exhibited. One is induced by the classical $sl(2)$ and the other by
the distinct one of ${\cal U}_{h}(sl(2))$. Automorphisms related to the two periods
of the elliptic functions involved are constructed. Translations of one 
generator by half and quarter periods lead to interesting results in this
context. Possibilities of applications are discussed briefly.
\end{abstract}

\newpage

\section{Introduction}
In [1] a nonlinear map relating the generators of the Jordanian ${\cal 
U}_{h}(sl(2))$ and $sl(2)$ was presented. Further applications of this 
mapping were studied in [2]. ( A fairly complete list of references  
concerning ${\cal U}_{h}(sl(2))$ can be found in [1] and [2].) Here a 
two-parametric generalization involving Jacobian elliptic functions of 
${\cal U}_{h}(sl(2))$ to ${\cal U}_{(h,k)}(sl(2))$ is presented. For
$k^2 =1$ one gets back ${\cal U}_{h}(sl(2))$. Setting now $h=0$ one 
gets $sl(2)$. The doubly periodic elliptic funtions involved lead to 
new interesting features.

Though two distinct induced Hopf structures will be exhibited for ${\cal U}_{(h,k)}
(sl(2))$, our main interest lies elsewhere. Exploring the enveloping algebra of
$sl(2)$ one encounters particularly interesting triplets leading to closed
nonlinear algebras with remarkable properties. (Here the term nonlinear
signifies that some of the commutators of the members of such a triplet give
nonlinear functions of them.) One may try to explore higher dimensional
algebras from this point of view. But so far our consructions have been
limited to $sl(2)$ and in that again mostly to two broad, complementary
classes. In one [3,4] the nonliearity arises through the generator
$J_0$ diagonalizable on the space of states spanning highest weight
irreducible representations of the basic triplet $(J_+,J_-,J_0)$. In the
other [1,2] this role is assumed by $J_+$ (or alternatively by $J_-$)
which is nilpotent on such spaces. These two types of nonlinearity can be
brought toghether. Such a consruction will be presented elsewhere [5],
combining the standard q-deformation  ${\cal U}_{q}(sl(2)$ with h-deformation
leading to ${\cal U}_{(q,h)}(sl(2))$.

Here we stay within the second class in generalizing  ${\cal U}_{h}(sl(2))$.
The nonlinear map is first presented and then the Hopf structures and
certain special automorphisms. A discussion of some interesting features 
and possibilities of applications will follow.

\section{The constuction as an invertible nonlinear map}

Let $(J_+,J_-,J_0)$ be the generators of $sl(2)$ satisfying
\begin{eqnarray}
&&[J_0,J_{\pm}]=\pm J_{\pm}\\
&&[J_+,J_-]= 2J_0
\end{eqnarray}
Define the set $({\hat X},{\hat Y},J_0)$ through the Jacobian elliptic 
functions [6,7] and the parameters $(h,k)$ as follows.
Let
\begin{eqnarray}
&&{h\over 2}J_+ = sn({h\over 2}{\hat X},k)
\end{eqnarray}
so that
\begin{eqnarray}
&&{h\over 2}{\hat X}= sn^{-1}({h\over 2}J_+,k).
\end{eqnarray}
Let
\begin{eqnarray}
&&{\hat Y}= g(J_+) J_-g(J_+)
\end{eqnarray}
where
\begin{eqnarray}
&&g(J_+) = \biggl(( 1 - ({h\over 2}J_+)^2 )( 1 - {k^2}({h\over 2}J_+)^2
)\biggr)^{1\over 4}\\
&&\phantom{g(J_+)}= \biggl( cn({h\over 2}{\hat X},k) dn({h\over 2}{\hat
X},k)\biggr)^{1\over 2}\\ && \phantom{g(J_+)}\equiv g(X)
\end{eqnarray}
so that
\begin{eqnarray}
&&J_- = g^{-1} (\hat X) {\hat Y} g^{-1}(\hat X).
\end{eqnarray} 

This map is to be understood here (at least to start with) in the context of
the standard highest weight representations of $sl(2)$ of dimention $(2j+1)$
for $j$ (half)integer. Since $J_+$ is nilpotent on such spaces the series 
developments of the elliptic functions involved are truncated. No convergence
problems arise. Thus, for example, from (3) and (4) setting $j ={ 5\over 2}$
keeping terms upto the fifth order only
\begin{eqnarray}
&&{h\over 2}J_+= {h\over 2}{\hat X} - {1\over {3!}}(1+ k^2){({h\over 2}{\hat
X})}^3 + {1\over {5!}}(1+14k+k^4){({h\over 2}{\hat X})}^5\\
&&{h\over 2}{\hat X}= {h\over 2}J_+ +{1\over {3!}}(1+k^2){({h\over 2}J_+)}^3
+{1\over {5!}}(9+6k^2+9k^4){({h\over 2}J_+)}^5
\end{eqnarray}
The series for other elliptic functions and their inverses are also truncated
analogously. From (1) to (9) one obtains for the triplet $({\hat X},{\hat
Y},J_0)$ the commutators
\begin{eqnarray}
&&[{\hat X},{\hat Y}]= 2J_0\\
&&[J_0,{\hat X}]= {2\over h}\biggl({sn({h\over 2}{\hat X},k)\over
cn({h\over 2}{\hat X},k)dn({h\over 2}{\hat X},k)}\biggr) \equiv G({\hat X})\\
&&[J_0,{\hat Y}]= -{1\over 2}(f({\hat X}){\hat Y}+{\hat Y}f({\hat X}))\\
&&f({\hat X})= {1 - k^2(sn({h\over 2}{\hat X},k))^4\over
(cn({h\over 2}{\hat X},k)dn({h\over 2}{\hat X},k))^2}\\
&&\phantom{f({\hat X})}={ 2\over sn(h{\hat X},k)}\biggl({sn({h\over 2}
{\hat X},k)\over
cn({h\over 2}{\hat X},k)dn({h\over 2}{\hat X},k)}\biggr)\\
&& \phantom{f({\hat X})}={ 1 - k^2({h\over 2}J_+)^4\over (1-({h\over
2}J_+)^2)( 1 - k^2({h\over 2}J_+)^2)}
\end{eqnarray}

For $k^2=1$, denoting
\begin{eqnarray}
&&({\hat X})_{k^2=1}= X\\
&&({\hat Y})_{k^2=1}= Y\\
&&{h\over 2}X= \hbox{arctanh}({h\over 2}J_+)\\
&&Y=\biggl (1-({h\over 2}J_+)^2\biggr)^{1\over 2}J_-\biggl(1-({h\over
2}J_+)^2\biggr)^{1\over 2} \end{eqnarray}
This is the nonlinear map introduced in [1] to obtain
${{\cal U}_h}(sl(2))$ satisfying 
\begin{eqnarray}
&&[X,Y]=2J_0\\
&&[J_0,X]= {1\over h}\sinh(hX)\\
&&[J_0,Y]= -{1\over 2}\biggl(\cosh(hX)Y+ Y\cosh(hX)\biggr)
\end{eqnarray} 
One can now express $({\hat X},{\hat Y})$ also through the map 
 \begin{eqnarray}
&&{h\over 2}{\hat X}= sn^{-1}(\tanh({h\over 2}X),k)\\
&&{\hat Y}= \biggl({1- k^2(\tanh({h\over2}X))^2\over 1- 
(\tanh({h\over 2}X))^2}\biggr)^{1\over 4}Y\biggl({1- k^2(\tanh({h\over2}X))^2
\over 1- 
(\tanh({h\over 2}X))^2}\biggr)^{1\over 4}
\end{eqnarray}
The unique Casimir for $sl(2)$ can be expressed in terms of the different
triplets as follows
\begin{eqnarray}
&&C= J_-J_+ + J_0(J_0+1)\\
&&\phantom{C}={2\over h}\cosh({h\over 2}X)Y \sinh({h\over 2}X)+J_0(J_0+1)\\
&&\phantom{C}={2\over h}(cn({h\over 2}{\hat X},k)dn({h\over 2}{\hat X},k)
)^{-{1\over 2}}{\hat  Y}(cn({h\over 2}{\hat X},k)dn({h\over 2}{\hat X},
k))^{-{1\over 2}}sn({h\over 2}{\hat X},k)+J_0(J_0+1)\qquad
\end{eqnarray}

The matrix elements of ${\hat X},{\hat Y}$ on the standard basis states
$\mid j,m \rangle$ of $sl(2)$ are obtained as follows. One has
\begin{eqnarray}
&&J_+\mid j,m \rangle = a_m\mid j,m +1\rangle\\
&&J_-\mid j,m \rangle = a_{m-1}\mid j,m-1 \rangle\\
&&J_0\mid j,m \rangle = m\mid j,m \rangle\\
&&a_m = ((j-m)(j+m+1))^{1\over2} ,\qquad (m = -j,\cdots,j)
\end{eqnarray}

Using the series developments [6,7] for the elliptic functions truncated by
\begin{eqnarray}
&&{J_+}^{(2j+1)}\mid j,m \rangle = 0
\end{eqnarray}
one obtains
\begin{eqnarray}
&&{\hat X}\mid j,m \rangle = a_m\mid j,m+1 \rangle + ({h\over2})^2{1\over 3!}
(1+k^2)(\prod_{i=0}^2 a_{m+i})\mid j,m+3 \rangle \nonumber \\
&& \phantom{ {\hat X}\mid j,m \rangle  }+({h\over2})^4{1\over 5!}
(1+14k^2+k^4)(\prod_{i=0}^4 a_{m+i})\mid j,m+5 \rangle + \cdots
 \end{eqnarray}

Similarly one can evaluate
\begin{eqnarray}
&&{\hat Y}\mid j,m \rangle = \biggl\{(1-{1\over 4}(1+k^2)({h\over 2}J_+)^2
-{1\over 32}(3+k^2)({h\over 2}J_+)^4\cdots)J_-\nonumber \\
&&\phantom{{\hat Y}\mid j,m \rangle }(1-{1\over
4}(1+k^2)({h\over 2}J_+)^2 -{1\over 32}(3+k^2)({h\over
2}J_+)^4\cdots)\Biggl\}\mid j,m \rangle \end{eqnarray}

For $j={5\over 2}$ for example, the terms exhibited explicitly in the series
suffice.
One has for each triplet
\begin{eqnarray}
&&C\mid j,m \rangle = j(j+1)\mid j,m \rangle .  
 \end{eqnarray} 

It can be shown that (12),(13) and (14) are compatible with the Jacobi identity.
 The crucial relation that assures this is $f({\hat X}) = {d\over d{\hat X}}
G({\hat X})$. 

\section{ The two induced Hopf sructures}

Let us now indicate the induced Hopf structures corresponding to the maps
(4), (5) and (25), (26) respectively. It is sufficient to consider the coproducts
as illustrations. The counits and the antipodes can be treated analogously.

For $sl(2)$ one has
\begin{eqnarray}
&&\Delta J_i = J_i \otimes 1+1 \otimes J_i, \qquad(i=\pm,0)
\end{eqnarray}
These induce, through (4), (5) and (6), the first sructure. Namely,
\begin{eqnarray}
&&\Delta_1{\hat X} = {2\over h}sn^{-1}({h\over 2}\Delta J_+,k)\\
&&\Delta_1{\hat Y} = g(\Delta J_+)(\Delta J_-)g(\Delta J_+)\\
&&\Delta J_0 = J_0 \otimes 1+1 \otimes J_0
\end{eqnarray}
Inverting the map (4), (5) one can express the righthand sides in terms of
$({\hat X},{\hat Y},J_0)$. Setting $k^2 = 1$ one obtains the induced coproducts
$(\Delta_1)$ of ${\cal U}_h(sl(2))$.

But as is wellknown, ${\cal U}_h(sl(2))$ has a distinct Hopf sructure leading to
a nontrivial $R$-mattrix. ( See [1] and the references cited there.) The
coproducts corresponding to this one are
\begin{eqnarray}
&&\Delta X = X \otimes 1 + 1 \otimes X\\
&&\Delta Y = Y \otimes e^{hX} + e^{-hX}  \otimes Y \\
&& \Delta J_0 = J_0 \otimes e^{hX} +e^{-hX}  \otimes J_0
\end{eqnarray}
Through (25), (26) one thus has a second set of induced coproducts 
  \begin{eqnarray}
&&\Delta_2 {\hat X}= {2 \over h}sn^{-1}(\tanh({h\over 2}\Delta X),k)\\
&&\Delta_2 {\hat Y}= \biggl({1- k^2(\tanh({h\over2}\Delta X))^2\over 1- 
(\tanh({h\over 2}\Delta X))^2}\biggr)^{1\over 4}\Delta Y\biggl({1-
k^2(\tanh({h\over2}\Delta X))^2 \over 1- 
(\tanh({h\over 2}\Delta X))^2}\biggr)^{1\over 4}\\
&&\Delta_2 J_0 = J_0 \otimes e^{hX} + e^{-hX}  \otimes J_0
 \end{eqnarray}

Inverting the map (25), (26) the righthand sides can be expressed again in terms
of $({\hat X},{\hat Y},J_0)$. Thus, for example
\begin{eqnarray}
&&\Delta_1{\hat X} = {2\over h}sn^{-1}\biggl(sn({h\over 2}{\hat X},k)\otimes 1
+ 1\otimes sn({h\over 2}{\hat X},k),k\biggr)\\
&&\Delta_2{\hat X} = {2\over h}sn^{-1}\biggl(\tanh(\tanh^{-1}sn({h\over 2}{\hat
X},k)\otimes 1 + 1\otimes \tanh^{-1}sn({h\over 2}{\hat X},k)),k\biggr)
\end{eqnarray}

One can similarly reexpress  $\Delta_1 {\hat Y}$ and $\Delta_2 {\hat Y}$.
For the non cocommutative $\Delta_2$ the R-matrices of ${{\cal U}_h}(sl(2)$ (of
which some examples can be found in [1] ) become relevant. But we will not
persue further these aspects in this paper. 

We have considered only induced Hopf structures. No new distinct one has 
been obtained.

\section{Automorphisms corresponding to half and quarter periods of
elliptic functions}

For comparison, we start by recapitulating the ${\cal U}_h(sl(2))$ automorphisms
 studied in [1]. Corresponding to the standard $sl(2)$ involution
\begin{eqnarray}
&&( J_+,J_-,J_0 ) \rightarrow ( -J_+,-J_-,J_0 )\\
&&( X,Y,J_0 ) \rightarrow ( -X,-Y,J_0 ).
\end{eqnarray}

Here the limit $h \rightarrow 0$ is straightforward. The situation is similar
for $({\hat X},{\hat Y},J_0)$ corresponding to (50).

But corresponding to half period of $\tanh$ one also has the automorphism
\begin{eqnarray}
&&( X,Y,J_0 ) \rightarrow ( X+{i\pi\over h},-Y,-J_0 )
\end{eqnarray}
which can be iterated. This has no straightforward $h \rightarrow 0$ 
(classical) limit. Since
\begin{eqnarray}
&&\tanh({h\over 2}( X + {i\pi\over h})) = \coth({h\over 2}X)
\end{eqnarray}
formally, from (20), (21) and (53) 
\begin{eqnarray} 
&&( J_+,J_-,J_0 ) \rightarrow ( J_+',J_-',J_0' )\\
&&J_+' = ({2\over h})^2J_+^{-1}\\
&&J_-' = ({h\over 2})^2J_+J_-J_+\\
&&J_0' = -J_0
\end{eqnarray}
This is an involution consistently with the fact that an iteration of (52)
\begin{eqnarray}
&&( X,Y,J_0 ) \rightarrow ( X+{{i2\pi}\over h},Y,J_0 )
\end{eqnarray}
is trivial at the level of $( J_+,J_-,J_0 )$.

The action of $J_+'$ on the $\mid j,m \rangle$ basis (see (30) to (33)) is no
longer well defined since ${a_m}^{-1}$ diverges for $m = j$. But the nonlinear
map assures that for $X$ there is only an additive imaginary term with
\begin{eqnarray}
&&X\mid j,j \rangle = {{i\pi}\over h}\mid j,j \rangle
\end{eqnarray}
This "half period automorphism" leads to a special class of $R$-matrices for
${\cal U}_h(sl(2))$ [1]. Let us now construct its generalization to the elliptic
case.

The two periods of the elliptic functions [6,7] are given in, terms of
$K$ and $K'$ where
\begin{eqnarray}
&&K = F({\pi\over 2},k) =  \int_0^{\pi\over 2}(I- k^2sin^2t)^{-{1\over2}}dt\\
&&K' =  F({\pi\over 2}, (1-k^2)^{1\over 2})
\end{eqnarray}
The primitive periods for $sn,cn,dn$ are respectively
\begin{eqnarray}
&&(4K,2iK') , \qquad (4K,2(K+iK')) ,\qquad  (2K,4iK')
\end{eqnarray}
 
Translations of $({h\over 2}{\hat X})$ by an entire common period leaves the
elliptic functions unchanged and hence the algebra is left evidently invariant.
But from the results (Table 7, page 350 [6]) corresponding to half and quarter
periods one can construct also the following automorphisms
\begin{eqnarray}
&&({\hat X},{\hat Y},J_0) \rightarrow ({\hat X}+{2\over h}iK',-{\hat Y},-J_0) \\
&&({\hat X},{\hat Y},J_0) \rightarrow ({\hat X}+{2\over h}(2K+iK'),-{\hat
Y},-J_0)
\end{eqnarray}

Note that the imaginary translation $iK'$ of the argument $({h\over 2}{\hat
X})$ corresponds to a quarter period for $dn$. But the combined effect of all
the factors leads to the follwing involutions (to be compared with (54)).
For nonzero values of $(h,k)$ one has with $\epsilon = \pm 1$ for (63), (64)
respectively
\begin{eqnarray}
&&(J_+,J_-,J_0) \rightarrow (\epsilon{1\over k}({2\over h})^2J_+^{-1},
\epsilon k({h\over 2})^2J_+J_-J_+,-J_0) .
\end{eqnarray}

Some relevant comments are added in the concluding discussions.

\section{Discussion}

We have presented our formalism, namely, the nonlinear map, the representations
and the automorphisms. Further developments can now be envisaged. 
The applications of $h$-deformation to $so(4)$ and $e(3)$ have been studied 
in [2]. One can try to generalize them in the context of $(h,k)$ deformation
 presented here.
 
One can try to select suitable classes of ansatz for Hamiltonians $H({\hat X},
{\hat Y},J_0)$ expressed directly as functions of these operators and explore
their properties. Our nonlinear map will provide a powerful tool in the study
of the spectra of such Hamiltonians.

In the context of magnetic fields perpendicular to a plane ${\cal U}_q(sl(2))$ has 
been used to construct Hamiltonians [8]. A classical limit has also been studied
[9]. In [10]$\;$a vector potential given by the Weierstrass zeta function has
been related to infinite dimentional representations of ${\cal U}_q(sl(2))$.

We have constructed two parametric deformations, here for $sl(2)$ and in [5]
in the context of ${\cal U}_q(sl(2))$. In our constructions the second parameter also
plays a full-fledged role. It appears in the algebra and is not relegated to 
the coproducts only. One can try to explore, not necessarily for magnetic
fields only, but in a broader context of models the potential role of these
parameters.

Let us finally come back to our automorphisms. In[11] outer automorphisms of 
infinite dimentional reprentations of $sl(2)$ involving $J_+^{-1}$ 
($s_+^{-1}$ in the notation of [11]) have been used to construct dynamical 
$r$-matrices. In (65) we have exhibited involutions leading to the inversion 
$J_+^{-1}$. At the level of the basic triplet $( J_+,J_-,J_0 )$ such inversions
 permit only infinite dimensional (no highest weight) representations. But for
 the triplet $({\hat X},{\hat Y},J_0)$ or $(X,Y,J_0)$ they still permit finite
 dimentional (though complex) representations. The explanation, here, is simple
enough. (When the eigenvalue of $ X$, say, on the highest weight state $ \mid j,j
\rangle$ is $({{i\pi}\over h})$ that of $\tanh({h\over 2}X)$ diverges.) But the
possibility and the consequences of having finite dimentional representations
 for suitably selected triplets in the enveloping algebra,when the realization
of the basic generators themselves permit only infinite dimentional ones,are
worth exploring more generally and systematically. Complementarily, going beyond 
the scope of this paper one can envisage independent and intrinsic explorations 
of infinite dimentional representations for the nonlinear algebras generated by 
$({\hat X},{\hat Y},J_0)$ or $(X,Y,J_0)$.

\smallskip

Elliptic functions appeared in deformations of $sl(2)$ in a quite different
fashion and with an extra generator in the works of Sklyanin [12]. References 
to recent "elliptic" deformations can be found in [4].

\newpage

It is a pleasure to thank Boucif Abdesselam for discussions and varied help.
 Our work related to $h$-deformations began in collaboration with and due 
to the initiative of R. Chakrabarti.

\end{document}